# PORTFOLIO CHOICE WITH JUMPS: A CLOSED-FORM SOLUTION

By Yacine Aït-Sahalia,[1] Julio Cacho-Diaz and T. R. Hurd[2]

*Princeton University, Princeton University and McMaster University*

We analyze the consumption-portfolio selection problem of an investor facing both Brownian and jump risks. We bring new tools, in the form of orthogonal decompositions, to bear on the problem in order to determine the optimal portfolio in closed form. We show that the optimal policy is for the investor to focus on controlling his exposure to the jump risk, while exploiting differences in the Brownian risk of the asset returns that lies in the orthogonal space.

**1. Introduction.** Economists have long been aware of the potential benefits of international diversification, while at the same time noting that the portfolios held by actual investors typically suffer from a home bias effect, meaning that those portfolios tend to be less diversified internationally than would be optimal according to portfolio choice theory [see, e.g., Solnik (1974) and Grauer and Hakansson (1987)]. One possible explanation is due to the risk of contagion across markets in times of crisis, which is a well-documented phenomenon [see, e.g., Longin and Solnik (2001), Ang and Chen (2002), Bae, Karolyi and Stulz (2003) and Hartmann, Straetmans and de Vries (2004)].

To generate this contagion, jumps of correlated sign across markets are a natural source of asymmetric correlation. Namely, when a downward jump occurs, negative returns tend to be experienced simultaneously across most markets, which then results in a high positive correlation in bear markets. When no jump occurs, the only source of correlation is that generated by the driving Brownian motions and will typically be much lower.

Studying dynamic portfolio choice has a long history, going back to Merton (1971). Subsequent papers have considered the portfolio problem, either

Received February 2008; revised June 2008.
[1]Supported in part by NSF Grants SBR-03-50772 and DMS-05-32370.
[2]Supported in part by NSERC and MITACS Canada.

in the one-period Markowitz setting or in the more complex dynamic Merton setting, when asset returns are generated by jump processes, including Poisson processes, stable processes or more general Lévy processes.[3] But so far, when jumps are included, the determination of an optimal portfolio has not been amenable to a closed-form solution, and this is a long-standing open problem in continuous-time finance. As a result, with $n$ assets, one must solve numerically an $n$-dimensional nonlinear equation. This is difficult, if not impossible, to do using existing methodologies: typical numerical examples consider the case $n = 1$ where a single risky asset is present. But with more efficient global markets, capital flows and a considerably larger number of potential assets to invest in, investors have more investment opportunities than ever before. In practice, we would certainly like to be able to solve models with hundreds if not thousands of assets.

This paper's main contribution is to make modeling choices and bring powerful yet simple mathematical tools, in the form of orthogonal decompositions, to bear on the problem: we will obtain the solution in closed form, thereby removing any restriction on the number of assets involved, as well as providing new insights and intuition into the structure of an optimal portfolio when jump risk is present.

---

[3]Early papers include Aase (1984), Jeanblanc-Picqué and Pontier (1990) and Shirakawa (1990). More recently, see Han and Rachev (2000) and Ortobelli et al. (2003) for a study of the Markowitz one-period mean-variance problem when asset returns follow a stable-Paretian distribution; Kallsen (2000) for a study of the continuous-time utility maximization in a market where risky security prices follow Lévy processes, and a solution (up to integration) for power, logarithmic and exponential utility using the duality or martingale approach; Choulli and Hurd (2001) give solutions up to constants of the primal and dual Merton portfolio optimization problem for the exponential, power and logarithmic utility functions when a risk-free asset and an exponential Lévy stock are the investment assets; Liu, Longstaff and Pan (2003) study the implications of jumps in both prices and volatility on investment strategies when a risk-free asset and a stochastic-volatility jump-diffusion stock are the available investment opportunities; Emmer and Klüppelberg (2004) study a continuous-time mean-variance problem with multiple assets; Madan (2004) derives the equilibrium prices in an economy with single period returns driven by exposure to explicit non-Gaussian systematic factors plus Gaussian idiosyncratic components. Cvitanić, Polimenis and Zapatero (2008) propose a model where the asset returns have higher moments due to jumps and study the sensitivity of the investment in the risky asset to the higher moments, as well as the resulting utility loss from ignoring the presence of higher moments. Das and Uppal (2004) evaluate the effect on portfolio choice of systemic risk, defined as the risk from infrequent events that are highly correlated across a large number of assets. They find that systemic risk reduces the gains from diversifying across a range of assets, and makes leveraged portfolios more susceptible to large losses. Ang and Bekaert (2002) consider a two-regime model in a discrete-time setting, one with low correlations and low volatilities, and one with higher correlations, higher volatilities and lower conditional means.

Specifically, we model asset returns as following exponential Lévy processes, which are a natural generalization of Merton's geometric Brownian motion model. With $n$ assets, the space of returns is $\mathbb{R}^n$. We show that the jump risk will occur in a well-defined portion of $\mathbb{R}^n$, say $\bar{V}$, and we then decompose $\mathbb{R}^n$ into $\bar{V}$ and the orthogonal space, $V^\perp$. Brownian risk, by contrast, occurs throughout $\mathbb{R}^n$. By adopting a factor structure on the Brownian variance-covariance of returns, we obtain a matching orthogonal decomposition of the Brownian risk which leads to a clean split of that Brownian risk into a portion occurring in $\bar{V}$, where it adds up to the jump risk, and a portion in $V^\perp$, where it is the sole source of risk. We are then able to distinguish between the optimal portfolio positions in the space $\bar{V}$ spanned by the jump risk (which the investor will attempt to limit) and those in the orthogonal space $V^\perp$ (where the investor will seek to exploit the opportunities arising from the traditional risk-return trade-off).

In our model, the structure of the Brownian variance-covariance matrix is taken to reflect the existence of one or more economic factors or sectors, each sector comprising a large number of related companies or countries. We start with the case where there is a single economic sector and jump, and then consider the more general case where the economy consists of $m$ sectors or regions of the world, each consisting of $k$ firms or countries. There, we allow for multiple jumps, potentially as many jump terms as there are sectors; those jumps can affect only their sector, or all or just some of the sectors, and do so to different degrees. The number of sectors can be as large as desired, providing a fair amount of flexibility in modeling assets with different characteristics and sources of jump risk.

Our closed-form solution provides new insights into the structure of the optimal portfolio in the presence of jumps. We show that the optimal investment policy can be summarized by a three-fund separation theorem, in the case of a single jump term. In general, if there exists enough cross-sectional variability in the expected excess returns, then the expected return and volatility of the portfolio value grow linearly in the number of assets. This happens because of exposure to the risky assets that the investor acquires in the space $V^\perp$. But the optimal policy in the subspace where jump risk lies, $\bar{V}$, is to control the portfolio's exposure to jumps by keeping it bounded as the number of assets grows. As a result, the overall exposure to jumps is dwarfed by the exposure to diffusive risk asymptotically in the number of assets $n$. Indeed, the additional investments in the risky assets are entirely in the direction $V^\perp$ that is orthogonal to the jump risk: they are all achieved with zero net additional exposure to the jump risk.

In other words, the optimal investment policy is to set the overall exposure to jump risk in $\bar{V}$ to the desired level, and then exploit, in the orthogonal space $V^\perp$, any perceived differences in expected returns and Brownian variances and covariances. But if the expected excess returns have little variability in the orthogonal space $V^\perp$, the opportunities for diversification effects are weak since controlling the exposure to jumps trumps other concerns (including the usual Brownian diversification policy).

The rest of the paper is organized as follows. In Section 2, we present our model of asset returns, and examine the investor's portfolio selection problem. In Section 3, we start for expository reasons with a simple case: a one-sector economy where the $n$ risky assets have the same jump size and introduce the mathematical tools we will use to derive the optimal portfolio weights in closed form. In Section 4, we provide a comparative statics analysis of the dependence of the optimal portfolio weights on the arrival intensity of the jumps, their magnitude, and the degree of risk aversion of the investor. In Section 5, we extend the model to an $m$-sector economy where sectors have their own sources of jump risk and solve the optimal portfolio problem in that case, again in closed form. In Section 6, we extend the model to other utility functions. Conclusions are in Section 7.

## 2. The portfolio selection model.

2.1. *Asset return dynamics.* Consider Merton's problem of maximizing the infinite-horizon expected utility of consumption by investing in a set of risky assets and a riskless asset. That is, the investor selects the amounts to be held in the $n$ risky assets and the riskless asset at times $t \in [0, \infty)$, as well as his consumption path. The available investment opportunities consist of a riskless asset with price $S_{0,t}$ and $n$ risky assets with prices $\mathbf{S}_t = [S_{1,t}, \ldots, S_{n,t}]'$.

For ease of exposition, we consider first the case of a single source of jump risk. In Section 5 below, we will generalize the model to include multiple jumps terms. Asset prices follow the exponential Lévy dynamics

$$(1) \qquad \frac{dS_{0,t}}{S_{0,t}} = r\, dt,$$

$$(2) \qquad \frac{dS_{i,t}}{S_{i,t-}} = (r + R_i)\, dt + \sum_{j=1}^{n} \sigma_{i,j}\, dW_{j,t} + J_i\, dY_t, \qquad i = 1, \ldots, n$$

with a constant rate of interest $r \geq 0$. $\mathbf{W}_t = [W_{1,t}, \ldots, W_{n,t}]'$ is an $n$-dimensional standard Brownian motion. $Y_t$ is a Lévy pure jump process with Lévy measure $\lambda \nu(dz)$, where $\lambda \geq 0$ is a fixed parameter and the measure $\nu$ satisfies $\int_{\mathbb{R}} \min(1, |z|) \nu(dz) < \infty$, so the jumps have finite variation. For any measurable subset $A$ of the real line, $\lambda \nu(A) = \lambda \int_A \nu(dz)$ measures the (possibly infinite) expected number of jumps, per unit of time, whose size belongs to $A$. If $\nu(A) = \infty$, it is necessarily because of small jumps, since $\nu(A) < \infty$ as long as $A$ does not contain 0. In other words, the expected number of jumps per unit of time of magnitude greater than any fixed $\varepsilon > 0$ is always finite. But there may be an infinite number of jumps of size less than $\varepsilon$.

The economy-wide jump amplitude $Y_t$ is scaled on an asset-by-asset basis by the scaling factor $J_i$. We assume that the support of $\nu$ is such that $S_i$ remains positive, consistent with the limited liability provision. We also assume that the jump process $Y$ and the individual Brownian motions in $\mathbf{W}$ are mutually independent.

In the special case where $Y$ is a compound Poisson process, there are a finite number of jumps of all sizes and $Y_t = \sum_{n=1}^{N_t} Z_n$, where $N_t$ is a scalar Poisson process with constant intensity parameter $\lambda > 0$, and the $Z_n$'s are i.i.d. random jump amplitudes with law $\nu(dz)$, that are independent of $N_t$. Then $\mathbf{S}$ follows a jump-diffusion.

In the absence of a jump term altogether, we obtain the usual geometric Brownian motion dynamics for $\mathbf{S}$.

The quantities $R_i$, $\sigma_{ij}$ and $J_i$ are constant parameters. We write $\mathbf{R} = [R_1, \ldots, R_n]'$, $\mathbf{J} = [J_1, \ldots, J_n]'$, and $\Sigma = \sigma \sigma'$, where

$$\sigma = \begin{pmatrix} \sigma_{1,1} & \cdots & \sigma_{1,n} \\ \vdots & \ddots & \vdots \\ \sigma_{n,1} & \cdots & \sigma_{n,n} \end{pmatrix}. \tag{3}$$

We assume that $\Sigma$ is a nonsingular matrix.

2.2. *Wealth dynamics and expected utility.* Let $\omega_{0,t}$ denote the percentage of wealth (or portfolio weight) invested at time $t$ in the riskless asset and $\omega_t = [\omega_{1,t}, \ldots, \omega_{n,t}]'$ denote the vector of portfolio weights in each of the $n$ risky assets, assumed to be adapted predictable càdlàg processes. The portfolio weights satisfy

$$\omega_{0,t} + \sum_{i=1}^{n} \omega_{i,t} = 1. \tag{4}$$

The investor consumes continuously at the rate $C_t$ at time $t$. In the absence of any income derived outside his investments in these assets, the investor's wealth, starting with the initial endowment $X_0$, follows the dynamics

$$\begin{aligned} dX_t &= -C_t \, dt + \omega_{0,t} X_t \frac{dS_{0,t}}{S_{0,t}} + \sum_{i=1}^{n} \omega_{i,t} X_t \frac{dS_{i,t}}{S_{i,t-}} \\ &= (rX_t + \omega_t' \mathbf{R} X_t - C_t) \, dt + X_t \omega_t' \sigma \, d\mathbf{W}_t + X_t \omega_t' \mathbf{J} \, dY_t. \end{aligned} \tag{5}$$

The investor's problem at time $t$ is then to pick the consumption and portfolio weight processes $\{C_s, \omega_s\}_{t \leq s \leq \infty}$ which maximize the infinite horizon, discounted at rate $\beta$, expected utility of consumption

$$V(X_t, t) = \max_{\{C_s, \omega_s; t \leq s \leq \infty\}} E_t \left[ \int_t^\infty e^{-\beta s} U(C_s) \, ds \right] \tag{6}$$

subject to the dynamics of his wealth (5), and with $X_t$ given. We will consider in detail in the rest of the paper the case where the investor has power utility, and briefly in Section 6 the cases of exponential and log utilities, respectively.

Using stochastic dynamic programming and the appropriate form of Itô's lemma for semimartingale processes, the Hamilton–Jacobi–Bellman equation characterizing the optimal solution to the investor's problem is

$$
\begin{aligned}
0 = \max_{\{C_t, \omega_t\}} \bigg\{ & e^{-\beta t} U(C_t) + \frac{\partial V(X_t, t)}{\partial t} \\
& + \frac{\partial V(X_t, t)}{\partial X}(rX_t + \omega_t' \mathbf{R} X_t - C_t) \\
& + \frac{1}{2} \frac{\partial^2 V(X_t, t)}{\partial X^2} X_t^2 \omega_t' \Sigma \omega_t \\
& + \lambda \int [V(X_t + X_t \omega_t' \mathbf{J} z, t) - V(X_t, t)] \nu(dz) \bigg\}
\end{aligned}
\tag{7}
$$

with the transversality condition $\lim_{t \to \infty} E[V(X_t, t)] = 0$ [see Merton (1969)].

The standard time-homogeneity argument for infinite-horizon problems gives that

$$
\begin{aligned}
e^{\beta t} V(X_t, t) &= \max_{\{C_s, \omega_s; t \leq s \leq \infty\}} E_t \bigg[ \int_t^\infty e^{-\beta(s-t)} U(C_s) \, ds \bigg] \\
&= \max_{\{C_{t+u}, \omega_{t+u}; 0 \leq u \leq \infty\}} E_t \bigg[ \int_0^\infty e^{-\beta u} U(C_{t+u}) \, du \bigg] \\
&= \max_{\{C_u, \omega_u; 0 \leq u \leq \infty\}} E_0 \bigg[ \int_0^\infty e^{-\beta u} U(C_u) \, du \bigg] \\
&\equiv L(X_t),
\end{aligned}
$$

which is independent of time. The third equality in this argument makes use of the fact that the optimal control is Markov. Thus $V(X_t, t) = e^{-\beta t} L(X_t)$ and (7) reduces to the following equation for the time-homogeneous value function $L$:

$$
\begin{aligned}
0 = \max_{\{C_t, \omega_t\}} \bigg\{ & U(C_t) - \beta L(X_t) + \frac{\partial L(X_t)}{\partial X}(rX_t + \omega_t' \mathbf{R} X_t - C_t) \\
& + \frac{1}{2} \frac{\partial^2 L(X_t)}{\partial X^2} X_t^2 \omega_t' \Sigma \omega_t \\
& + \lambda \int [L(X_t + X_t \omega_t' \mathbf{J} z) - L(X_t)] \nu(dz) \bigg\}
\end{aligned}
\tag{8}
$$

with the transversality condition

$$
\lim_{t \to \infty} E[e^{-\beta t} L(X_t)] = 0.
\tag{9}
$$

The maximization problem in (8) separates into one for $C_t$, with first-order condition

$$\frac{\partial U(C_t)}{\partial C} = \frac{\partial L(X_t)}{\partial X}$$

and one for $\omega_t$:

(10)
$$\max_{\{\omega_t\}}\left\{\frac{\partial L(X_t)}{\partial X}\omega_t'\mathbf{R}X_t + \frac{1}{2}\frac{\partial^2 L(X_t)}{\partial X^2}X_t^2\omega_t'\Sigma\omega_t \right.$$
$$\left. + \lambda\int[L(X_t + X_t\omega_t'\mathbf{J}z) - L(X_t)]\nu(dz)\right\}.$$

Given wealth $X_t$, the optimal consumption choice is therefore

(11)
$$C_t^* = \left[\frac{\partial U}{\partial C}\right]^{-1}\left(\frac{\partial L(X_t)}{\partial X}\right).$$

In order to determine the optimal portfolio weights, wealth and value function, we need to be more specific about the utility function $U$.

2.3. *Power utility.* Consider an investor with power utility, $U(c) = c^{1-\gamma}/(1-\gamma)$ for $c > 0$ and $U(c) = -\infty$ for $c \leq 0$ with CRRA coefficient $\gamma \in (0,1) \cup (1,\infty)$. (In Section 6, we briefly treat the exponential and log utility cases.) We will look for a solution to (8) in the form

(12)
$$L(x) = K^{-\gamma}x^{1-\gamma}/(1-\gamma),$$

where $K$ is a constant, so that

(13)
$$\frac{\partial L(x)}{\partial x} = (1-\gamma)L(x)/x, \qquad \frac{\partial^2 L(x)}{\partial x^2} = -\gamma(1-\gamma)L(x)/x^2.$$

Then (8) reduces to

(14)
$$0 = \max_{\{C_t,\omega_t\}}\Big(U(C_t) - \beta L(X_t) + (1-\gamma)L(X_t)(\omega_t'\mathbf{R} + r)$$
$$- (1-\gamma)C_t\frac{L(X_t)}{X_t}$$
$$- \frac{1}{2}\gamma(1-\gamma)L(X_t)\omega_t'\Sigma\omega_t$$
$$+ \lambda\int[(1+\omega_t'\mathbf{J}z)^{1-\gamma}L(X_t) - L(X_t)]\nu(dz)\Big),$$

that is,

(15)
$$0 = \min_{\{C_t,\omega_t\}}\Big(-\frac{U(C_t)}{(1-\gamma)L(X_t)} + \frac{\beta}{(1-\gamma)} - (r+\omega_t'\mathbf{R}) + \frac{C_t}{X_t}$$
$$+ \frac{1}{2}\gamma\omega_t'\Sigma\omega_t - \frac{\lambda}{(1-\gamma)}\int[(1+\omega_t'\mathbf{J}z)^{1-\gamma} - 1]\nu(dz)\Big),$$

after division by $-(1-\gamma)L(X_t) < 0$, so that max becomes min.

2.4. *Optimal policies.* The optimal policy for the portfolio weights $\omega_t$ is

$$\omega_t^* = \underset{\{\omega_t\}}{\arg\min}\, g(\omega_t), \tag{16}$$

where the functions

$$g(\omega) = -\omega'\mathbf{R} + \frac{\gamma}{2}\omega'\Sigma\omega + \lambda\psi(\omega'\mathbf{J}) \tag{17}$$

and

$$\psi(\omega'\mathbf{J}) = -\frac{1}{(1-\gamma)}\int [(1+\omega'\mathbf{J}z)^{1-\gamma} - 1]\nu(dz) \tag{18}$$

are both convex.

Since $(\gamma, \mathbf{R}, \Sigma, \mathbf{J})$ are constant, the objective function $g$ is time independent, so it is clear that any optimal solution will be time independent. Furthermore, the objective function is state independent, so any optimal solution will also be state independent. In other words, any optimal $\omega_t^*$ will be a constant $\omega^*$ independent of time and state. Finally, the objective function $g$ is strictly convex, goes to $+\infty$ in all directions, and hence always has a unique minimizer. In the pure diffusive case, $\lambda = 0$ and we obtain of course the familiar Merton solution

$$\omega^* = \frac{1}{\gamma}\Sigma^{-1}\mathbf{R}. \tag{19}$$

As to the optimal consumption policy, with $[\partial U/\partial C]^{-1}(y) = y^{-1/\gamma}$ and $\partial L(x)/\partial x = K^{-\gamma}x^{-\gamma}$ in (11), we obtain

$$C_t^* = KX_t \tag{20}$$

provided that $X_t > 0$. Next, we evaluate (15) at the optimal policies $(C_t^*, \omega^*)$ to identify the constant $K$:

$$K = \frac{\beta}{\gamma} - \frac{(1-\gamma)}{\gamma}[\omega^{*\prime}\mathbf{R} + r] + \frac{1}{2}(1-\gamma)\omega^{*\prime}\Sigma\omega^* + \frac{(1-\gamma)\lambda}{\gamma}\psi(\omega^{*\prime}\mathbf{J}). \tag{21}$$

The constant $K$ will be fully determined once we have solved below for the optimal portfolio weights, $\omega^*$.

Finally, we have to check that the transversality condition is satisfied. By plugging the optimizers $X^*$ and $C_t^*$ into (6), and then taking expectations, one finds

$$E[V(X_t^*, t)] = E\left[\int_t^\infty e^{-\beta s} U(C_s^*)\, ds\right].$$

Now $e^{-\beta s}U(C_s^*) = KV(X_s^*, s)$ from which it follows that $E[V(X_t^*, t)]$ solves $df/dt = -Kf$ and hence decays exponentially to zero as $t \to \infty$, for any $K > 0$.

By plugging the constants $K, \omega^*$ into (5), we see that an investor with power utility who selects this optimal portfolio will achieve a wealth process $X_t^*$ which follows a geometric Lévy process.

### 3. Optimal portfolio in a one-sector economy with homogenous jumps.

3.1. *The orthogonal decomposition.* To begin, we consider the simplest possible case, where the $n$ risky assets have identical jump size characteristics

$$(22) \qquad \mathbf{J} = \bar{J}\mathbf{1},$$

where $\bar{J}$ is a scalar, and $\mathbf{1}$ is the $n$-vector $\mathbf{1} = [1, \ldots, 1]'$. The key to finding the solution in closed form is to pick a basis of the space of returns, $\mathbb{R}^n$, that isolates the direction of the jump vector. In this simple case, jumps are parallel to the vector $\mathbf{1}$, so we are led to consider the orthogonal decomposition $\mathbb{R}^n = \bar{V} \oplus V^\perp$ where $\bar{V}$ is the span of $\mathbf{1}$ and $V^\perp$ is the orthogonal hyperplane.

By construction, all jump risk is contained in $\bar{V}$. However, Brownian risk occurs throughout $\mathbb{R}^n$. We will then have to figure out which part of the Brownian risk takes place in $\bar{V}$ (call it $\bar{\Sigma}$) and which part takes place in $V^\perp$ (call it $\Sigma^\perp$). Indeed, the key to going further, and in particular to finding solutions in closed form, is to make a further assumption about the relation between the Brownian risk and the subspace $\bar{V}$. We therefore make the following assumption on $\Sigma$:

ASSUMPTION 1 (*Orthogonal decomposition assumption*). We assume that $\bar{V}$ and $V^\perp$ are invariant subspaces of the matrix $\Sigma$. In other words, we assume that $\Sigma = \bar{\Sigma} + \Sigma^\perp$ maps $\bar{V}$ onto $\bar{V}$ and $V^\perp$ onto $V^\perp$.

This is where assuming a factor structure for the Brownian variance-covariance matrix of returns $\Sigma$ will be useful, as it will make the decomposition $\Sigma = \bar{\Sigma} + \Sigma^\perp$ explicit, as we will see below. Throughout the paper, including the multisector case of Section 5 where the jump and factor structures are much richer, we will always decompose the returns space $\mathbb{R}^n$ into the direction of the jump vector(s), where both Brownian and jump risks occur (for now $\bar{V}$ is simply $\mathbf{1}$) and the orthogonal space ($V^\perp$, where only Brownian risk occurs). This decomposition is illustrated in Figure 1 in the simple case we are starting with.

As will become clear below, this decomposition leads to a separation of the problem into two subproblems: one in $\bar{V}$, and one in $V^\perp$. Because $V^\perp$ contains only Brownian risk, the solution in that space will be the Merton no-jump solution, except that it will be based not on the full $\Sigma$, but just on the fraction $\Sigma^\perp$ of Brownian risk that lies in that space. That part of the solution will be explicit, as the Merton no-jump solution always is. In $\bar{V}$, we will have to deal with both sources of risk.

But we will already have achieved dimension-reduction: the solution for the portfolio weights in the large-dimensional space, $V^\perp$, is known, and the remaining unknowns in $\bar{V}$ are low-dimensional (one-dimensional here in fact). With more structure on the dynamics generating the asset returns, namely specific assumptions on the jump measure $\nu$, we will show how that last unknown (the optimal portfolio policy in $\bar{V}$) can also be obtained in closed form.

So, the key is to work with the orthogonal decomposition of $\mathbb{R}^n$ that is suggested by the jumps. With that in mind, we therefore decompose the vector of excess returns as follows:

$$\mathbf{R} = \bar{R}\mathbf{1} + \mathbf{R}^\perp, \tag{23}$$

where $\bar{R}$ is a scalar and $\mathbf{R}^\perp$ is an $n$-vector orthogonal to $\mathbf{1}$.

As for $\Sigma$, we assume for now the one-factor structure

$$\Sigma = v^2 \begin{pmatrix} 1 & \rho & \cdots \\ \rho & \ddots & \rho \\ \cdots & \rho & 1 \end{pmatrix}, \tag{24}$$

where $v^2 > 0$ is the variance of the returns generated by the diffusive risk, and $-1/(n-1) < \rho < 1$ is their common correlation coefficient.

The $\Sigma$ matrix decomposes on $\mathbb{R}^n = \bar{V} \oplus V^\perp$ as follows:

$$\Sigma = \underbrace{\kappa_1 \frac{1}{n}\mathbf{11}'}_{=\bar{\Sigma}} + \underbrace{\kappa_2 \left(\mathbf{I} - \frac{1}{n}\mathbf{11}'\right)}_{=\Sigma^\perp}, \tag{25}$$

where $\mathbf{I}$ denotes the $n \times n$ identity matrix and

$$\kappa_1 = v^2 + v^2(n-1)\rho, \tag{26}$$

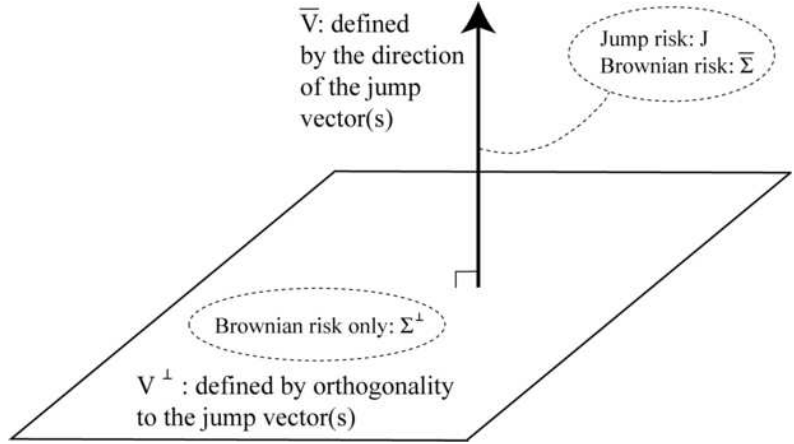

FIG. 1. *Orthogonal decomposition of the returns space $\mathbb{R}^n = \bar{V} \oplus V^\perp$ where $\bar{V}$ is the span of the jump vector(s) and $V^\perp$ is the orthogonal hyperplane. By construction, only $\bar{V}$ contains jump risk. However, Brownian risk is contained throughout $\mathbb{R}^n$: $\bar{V}$ contains the part $\bar{\Sigma}$, and $V^\perp$ the part $\Sigma^\perp$.*

$$\kappa_2 = v^2(1-\rho) \tag{27}$$

are the two distinct eigenvalues of $\Sigma$, $\kappa_1$ with multiplicity 1 and eigenvector **1** and $\kappa_2$ with multiplicity $n-1$.

We then search for the optimal portfolio vector $\omega$ in the form

$$\omega = \bar{\omega}\mathbf{1} + \omega^{\perp}, \tag{28}$$

where $\bar{\omega}$ is scalar and $\omega^{\perp}$ is an $n$-vector orthogonal to **1**.

3.2. *Optimal portfolio separation.* Given the decomposition (25), we see from (16)–(17) that the optimal $\bar{\omega}^*$ and $\omega^{\perp *}$ must satisfy

$$(\omega^{\perp *}, \bar{\omega}^*) = \arg\min_{\{\omega^{\perp}, \bar{\omega}\}} \left\{ -n\bar{\omega}\bar{R} + \frac{1}{2}\gamma n\bar{\omega}^2\kappa_1 + \lambda\psi(n\bar{\omega}\bar{J}) \right. \tag{29}$$
$$\left. - \omega^{\perp\prime}\mathbf{R}^{\perp} + \frac{1}{2}\gamma\omega^{\perp\prime}\kappa_2\left(\mathbf{I} - \frac{1}{n}\mathbf{11}'\right)\omega^{\perp} \right\}.$$

And we now see that this separates into two distinct optimization problems, one for $\omega^{\perp}$ and one for the scalar $\bar{\omega}$:

$$\begin{cases} \omega^{\perp *} = \arg\min_{\omega^{\perp}} \{g^{\perp}(\omega^{\perp})\}, \\ \bar{\omega}^* = \arg\min_{\bar{\omega}} \{\bar{g}(\bar{\omega})\}, \end{cases} \tag{30}$$

where

(31) $g^{\perp}(\omega^{\perp}) = - \underbrace{\omega^{\perp\prime}\mathbf{R}^{\perp}}_{\text{Part of return in } V^{\perp}} + \underbrace{\frac{1}{2}\gamma\kappa_2\omega^{\perp\prime}\omega^{\perp}}_{\text{Part of Brownian risk in } V^{\perp}},$

(32) $\bar{g}(\bar{\omega}) = - \underbrace{n\bar{\omega}\bar{R}}_{\text{Part of return in } \bar{V}} + \underbrace{\frac{1}{2}\gamma n\bar{\omega}^2\kappa_1}_{\text{Part of Brownian risk in } \bar{V}} + \underbrace{\lambda\psi(n\bar{\omega}\bar{J})}_{\text{Jump risk}}.$

In $V^{\perp}$, things are simple because there is no jump risk by construction. The first-order condition for minimizing (31) is

$$-\mathbf{R}^{\perp} + \gamma\kappa_2\omega^{\perp *} = 0,$$

whose solution is

$$\omega^{\perp *} = \frac{1}{\gamma\kappa_2}\mathbf{R}^{\perp} = \frac{1}{\gamma v^2(1-\rho)}\mathbf{R}^{\perp}. \tag{33}$$

This is nothing else than the Merton no-jump solution, but *in restriction* to the space $V^{\perp}$. In other words, $\omega^{\perp *}$ is the solution in the no-jump space $V^{\perp}$ for the Merton problem with Brownian risk $\Sigma^{\perp}$ and the vector of excess returns $\mathbf{R}^{\perp}$. The investor is going after excess returns, subject to the usual continuous-risk provisions.

To solve (30) for the last remaining unknown $\bar{\omega}^*$, we have to contend with both Brownian and jump risks. Below, we will show how to determine $\bar{\omega}^*$ in closed form under further assumptions on the distribution of the jumps and the investor's utility function. But for now, we see that the optimal portfolio choice is characterized by

$$\omega^* = \bar{\omega}^* \mathbf{1} + \frac{1}{\gamma v^2 (1-\rho)} \mathbf{R}^\perp. \tag{34}$$

3.3. *Large $n$ asymptotics.* As far as the optimal solution for $\bar{\omega}$ is concerned, with the change of variable $\varpi_n = n\bar{\omega}$, we see that

$$\varpi_n^* = \arg\min_{\varpi_n} \left\{ -\varpi_n \bar{R} + \frac{1}{2}\gamma \varpi_n^2 \kappa_1/n + \lambda \psi(\varpi_n \bar{J}) \right\}. \tag{35}$$

Letting $n \to \infty$, we have that $\kappa_1/n \to v^2 \rho$ and so $\varpi_n^*$ converges to a finite constant $\varpi_\infty^*$ given by

$$\varpi_\infty^* = \arg\min_{\varpi_\infty} \left\{ -\varpi_\infty \bar{R} + \frac{1}{2}\gamma \varpi_\infty^2 v^2 \rho + \lambda \psi(\varpi_\infty \bar{J}) \right\}. \tag{36}$$

This convergence is illustrated in Figure 3.

An investor who selects this optimal portfolio will achieve a wealth process $X_t^*$ which follows a geometric Lévy process with characteristic triple $(b, c, f)$ such that the drift $b$ and volatility $c$ are $O(n)$, as long as there is sufficient cross-sectional dispersion of excess returns, in the sense that $\|\mathbf{R}^\perp\|^2 = \mathbf{R}^{\perp\prime} \mathbf{R}^\perp = O(n)$, while the Lévy jump measure $f$ remains $O(1)$ as $n \to \infty$.

This means that expected excess returns $\mathbf{R}^\perp$ lead the optimal portfolio value to increasing expected returns $b$ and variance $c$, both growing linearly in the number of assets. On the other hand, as $n$ grows, the exposure to contagion jumps remains bounded, and is dwarfed by the exposure to diffusive risk. Indeed, the investment in the risky assets due to the expected return compensation $\mathbf{R}^\perp$ is entirely in the direction of $\omega^\perp$, which is orthogonal to $\mathbf{J}$. Since all the assets have the same response to the jump risk under the simple model (22), this is all done in a long-short manner: the sum of the wealth percentages invested in $V^\perp$ is 0. So the additional amounts invested in the risky assets are all achieved without increasing the net exposure to the jump risk. In fact, since $\varpi_n^* \to \varpi_\infty^*$ finite, we have asymptotically $\bar{\omega}^* \sim \varpi_\infty^*/n$ so the investor exploits an increase in the number of assets to reduce, ceteris paribus, his net exposure to the risky assets.

3.4. *A three-fund separation theorem.* It is clear from the analysis above that there exist three funds such that an investor should be indifferent between the original $n$ assets (plus the riskless asset) or just the three funds;

the proportions of each fund invested in the individual assets depend only on the asset return characteristics and not on investor preferences; and the investor's optimal demands for the three funds depend on his preferences. By contrast, both in the static one-period Markowitz model and in the dynamic Merton model, one needs two funds (the market portfolio and the riskless asset) to replicate the optimal asset allocation.

In our simple one-sector model with homogenous jump sizes, the three funds are:

- A long-short fund LS which holds proportions $\delta_1 = \frac{1}{v^2(1-\rho)}\mathbf{R}^\perp$ of the $n$ risky assets. This fund is long-short because it contains no net exposure to the risky assets, since the sum of the weights invested in the various assets by fund LS is zero. Said differently, any investment in fund LS generates no new exposure to jump risk.
- An equally weighted fund EQ holding a portfolio $\delta_2 = \frac{1}{n}\mathbf{1}$ of the $n$ risky assets. Investment in this fund generates net exposure to jump risk. How much of that risk the investor is willing to take is controlled by the quantity $\bar{\omega}^*$.
- A risk-free fund RF holding the riskless asset only.

Funds LS and RF provide an efficient diffusion risk-return optimization (for the part of the Brownian risk orthogonal to $\mathbf{1}$, that is, the part in $V^\perp$), while fund EQ allows investors to hedge against the jump risk, plus the part of the Brownian risk in the direction of $\mathbf{1}$, that is, the part in $\bar{V}$. Replicating the investor's optimal positions cannot be achieved by lumping together LS and EQ into a single market portfolio.

The fact that the "jump-hedging" fund is perfectly equally weighted (or equivalently the fact that the orthogonal fund is perfectly long-short) is of course due to the fact that in the model so far all assets are symmetrically exposed to the jump risk. The fact that we need three funds is due to the fact that so far we have only one source of jump risk. We will generalize all of this below in Section 5.

3.5. *Fully explicit portfolio weights.* Some special cases lead to a closed-form solution for the last remaining constant, $\varpi_n^*$, yielding fully closed-form solutions for the optimal portfolio weights. For this, we need to specify the utility function and the Lévy measure $\nu(dz)$ driving the common jumps; then we can compute the integral in (36). We provide a few examples.

Consider the case of a power utility investor with CRRA coefficient $\gamma = 2$, the Lévy measure satisfying a power law,

(37) $$\lambda\nu(dz) = \begin{cases} \lambda_+ \, dz/z, & \text{if } z \in (0,1], \\ -\lambda_- \, dz/z, & \text{if } z \in [-1,0), \end{cases}$$

and $\bar{J} \in (-1, 1)$. We define $\lambda_+ > 0$ as the intensity of positive jumps and $\lambda_- > 0$ as the intensity of negative jumps. Equation (36) specializes to

$$\varpi_\infty^* = \arg\min_\varpi f_\infty(\varpi), \tag{38}$$

where

$$f_\infty(\varpi) = -\varpi\bar{R} + \varpi^2 v^2 \rho + \lambda_+ \int_0^1 [(1 + \varpi\bar{J}z)^{-1} - 1] \, dz/z \tag{39}$$

$$- \lambda_- \int_{-1}^0 [(1 + \varpi\bar{J}z)^{-1} - 1] \, dz/z \tag{40}$$

$$= -\varpi\bar{R} + \varpi^2 v^2 \rho - \lambda_+ \log(1 + \varpi\bar{J}) - \lambda_- \log(1 - \varpi\bar{J}).$$

The first-order condition (FOC) for $\varpi$ is given by

$$-\bar{R} + 2\varpi v^2 \rho - \lambda_+ \bar{J}(1 + \varpi\bar{J})^{-1} + \lambda_- \bar{J}(1 - \varpi\bar{J})^{-1} = 0. \tag{41}$$

The optimal solution must satisfy

$$|\bar{J}\varpi_\infty^*| < 1; \tag{42}$$

otherwise, there is a positive probability of wealth $X_t$ becoming negative, which is inadmissible in the power utility case. The asymptotic solution to (41) is the unique root $\varpi_\infty^*$ satisfying the solvency constraint (42), and that solution is given by

$$\varpi_\infty^* = \frac{\bar{R}}{6v^2\rho}$$

$$+ \sqrt{\frac{4}{3}\left(\frac{\bar{R}^2}{12(v^2\rho)^2} + \frac{\lambda_+ + \lambda_-}{2v^2\rho} + \frac{1}{\bar{J}^2}\right)}$$

$$\times \cos\left(\frac{1}{3}\arccos\left(\left(\frac{\bar{R}}{6v^2\rho}\left(\frac{\bar{R}^2}{36(v^2\rho)^2} + \frac{\lambda_+ + \lambda_-}{4v^2\rho}\right.\right.\right.$$

$$\left.\left.- \frac{3(\lambda_+ - \lambda_-)}{2\bar{J}\bar{R}} - \frac{1}{\bar{J}^2}\right)\right)$$

$$\times \left(\left(\sqrt{\frac{1}{3}\left(\frac{\bar{R}^2}{12(v^2\rho)^2} + \frac{\lambda_+ + \lambda_-}{2v^2\rho} + \frac{1}{\bar{J}^2}\right)}\right)^3\right)^{-1}\right)$$

$$\left. + \frac{4}{3}\pi\right). \tag{43}$$

3.5.1. *Negative jumps.* Let us now consider the case $\bar{J} < 0$, $\lambda_- = 0$ and $\lambda = \lambda_+ > 0$ in order to capture the downward risk inherent in the types of jumps we are concerned about. Figure 2 plots $\varpi_\infty^*$ as a function of $\bar{J}$ and $\lambda$. Visual inspection of (43) reveals that $\bar{J}$ and $\lambda$ do not have a symmetric effect on the optimal portfolio weights. In the exact small-sample case, the optimal solution (43) becomes

$$(44) \qquad \varpi_n^* = \frac{2\kappa_1/|\bar{J}| + n\bar{R} - \sqrt{(2\kappa_1/|\bar{J}| - n\bar{R})^2 + 8n\kappa_1\lambda}}{4\kappa_1}.$$

Figure 3 plots the objective function, $f_n(\varpi) = \bar{g}(\varpi/n)$, and shows its convergence to $f_\infty(\varpi)$ as $n \to \infty$, along with $\arg\min f_n(\varpi) = \varpi_n^*$, converging to $\varpi_\infty^*$. It is also worth noting that

$$(45) \qquad \varpi_n^* < \frac{2\kappa_1/|\bar{J}| + n\bar{R} - |2\kappa_1/|\bar{J}| - n\bar{R}|}{4\kappa_1} = \min\left(\frac{n\bar{R}}{2\kappa_1}, \frac{1}{|\bar{J}|}\right),$$

so that the optimal investment in the risky assets is always less than what it would be in the absence of jumps, $n\bar{R}/(2\kappa_1)$. This is natural since $\bar{J} < 0$.

3.5.2. *Other examples.* Other cases that lead to a closed-form solution for $\varpi_n^*$ include power utility with $\gamma = 3$ and either power law jumps $\nu(dz) = \lambda dz/z$, or uniform jumps $\nu(dz) = \lambda dz$. In either case, the FOC is then a cubic equation, as (41), solvable in closed form using standard methods. Another case is one where the investor has log utility with jumps of a fixed size, $\nu(dz) = \delta(z = \bar{z})dz$, for some $\bar{z} \in [-1, 1]$. Then the FOC leads to a quadratic equation.

**4. Comparative statics.** Using the explicit solutions above, we can now investigate how the optimal portfolio responds to different jump intensities, jump sizes and degrees of risk aversion. This analysis is decidedly partial-equilibrium, in that we do not allow for changes in the assets' expected returns when other parameters of the model change. To save space, we consider only the negative jump example of Section 3.5.1, where $\bar{J} < 0$, $\lambda_- = 0$ and $n \to \infty$.

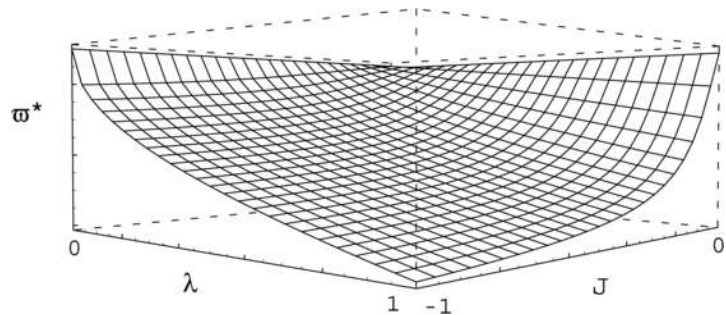

FIG. 2. *Optimal portfolio weight $\varpi_\infty^*$ as a function of $\bar{J}$ and $\lambda$.*

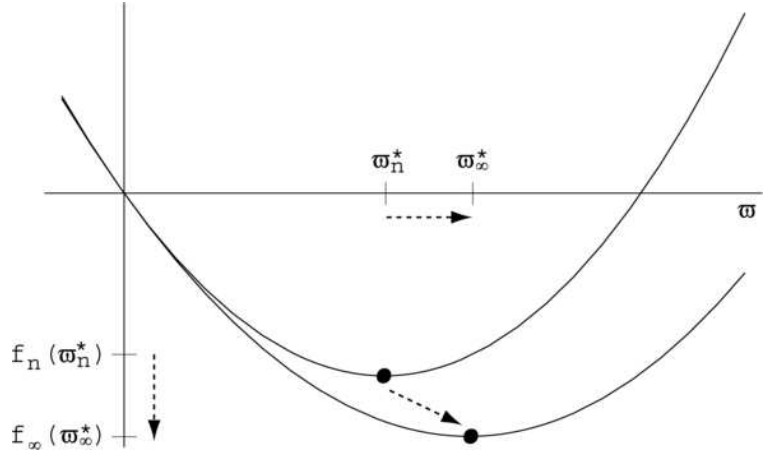

FIG. 3. *Scalar objective function used to determine the optimal portfolio weight $\varpi_n^*$ and its large asset asymptotic limit, $\varpi_\infty^*$.*

4.1. *Response to jumps of different arrival intensity.* We have

(46) $$\varpi_\infty^* \to -\infty \qquad \text{as } \lambda \to \infty,$$

(47) $$\varpi_\infty^* \to \min\left(\frac{\bar{R}}{2v^2\rho}, \frac{1}{|\bar{J}|}\right) \qquad \text{as } \lambda \to 0.$$

The first limit means that the investor will go short to an unbounded extent on all the risky assets if the arrival rate of the jumps goes to infinity. This is to be expected, since $\bar{J} < 0$ and we impose no short sale constraints. Further, $\varpi_\infty^*$ tends to $-\infty$ when $\lambda \to \infty$ at the following rate:

$$\varpi_\infty^* = -\frac{\sqrt{\lambda}}{\sqrt{2v^2\rho}}(1 + o(1)).$$

If, on the other hand, the jumps become less and less frequent, $\lambda \to 0$, then $\varpi_\infty^*$ tends to a finite limit driven by the diffusive characteristics of the assets, $\bar{R}/(2v^2\rho)$, which is the limit in the no-jump case, *unless* the jump size $\bar{J}$ is so large that the solvency constraint binds. This gives rise to a kink in the demand function.

If the solvency constraint is not binding, then the higher the variance of the assets and/or the more heavily correlated they are, the smaller the investment in each one of them. And the higher the expected excess return of the assets $\bar{R}$, the higher the amount invested. For a small perceived jump risk ($\lambda$ small), the optimal solution behaves like

$$\varpi_\infty^* = \min\left(\frac{\bar{R}}{2v^2\rho}, \frac{1}{|\bar{J}|}\right) + \frac{\bar{J}\lambda}{|\bar{R}+\bar{J}\lambda|} + o(\lambda).$$

4.1.1. *Jumps versus expected return trade-off.* The weights $\varpi_\infty^*$ are monotonic in $\lambda$, with

$$\frac{\partial \varpi_\infty^*}{\partial \lambda} = \frac{\bar{J}}{\sqrt{(2v^2\rho - \bar{J}\bar{R})^2 + 8\bar{J}(\bar{R}+\bar{J}\lambda)v^2\rho}} < 0.$$

If $\bar{R} > 0$, there exists a critical value $\tilde{\lambda}$ such that

(48) $\qquad\qquad\qquad \varpi_\infty^* > 0 \qquad \text{for } \lambda < \tilde{\lambda},$

(49) $\qquad\qquad\qquad \varpi_\infty^* \leq 0 \qquad \text{for } \lambda \geq \tilde{\lambda}.$

That is, as long as jumps do not occur too frequently ($\lambda < \tilde{\lambda}$), the investor will go long on the assets in order to capture their expected return, even though that involves taking on the (negative) risk of the jumps. When the jumps occur frequently enough ($\lambda \geq \tilde{\lambda}$), then the investor decides to forgo the expected return of the assets and focuses on canceling his exposure to the jump risk by going short these assets.

The critical value $\tilde{\lambda}$ takes a particularly simple form. It is given by

(50) $$\tilde{\lambda} = \frac{\bar{R}}{|\bar{J}| \int_0^1 z\nu(dz)}.$$

Clearly, the higher $\bar{R}$ relative to $|\bar{J}|$, the higher $\tilde{\lambda}$. And the smaller the expected value of $Z$, the bigger $\tilde{\lambda}$. With $\nu(dz) = dz/z$, we get

(51) $$\tilde{\lambda} = \frac{\bar{R}}{|\bar{J}|}.$$

This expression can be interpreted as an analogue to the Sharpe ratio, but for jump risk: excess return in the direction of jumps, $\bar{R}$, divided by a measure of the jump magnitude, $|\bar{J}|$.

Now, if $\bar{R} \leq 0$, then $\varpi_\infty^* \leq 0$ for every $\lambda \geq 0$. In that case, there is no point in ever going long those assets since both the expected return and the jump components negatively impact the investor's rate of return.

4.1.2. *Flight to quality.* The solution above can capture a well-documented empirical phenomenon. Imagine an investor currently in a normal, low jump risk environment, who receives information suggesting that the jumps become more likely. Starting from a situation where $\lambda < \tilde{\lambda}$, if the perception of the jump risk increases ($\lambda \uparrow \tilde{\lambda}$), then the optimal policy for the investor is to flee-to-quality, by reducing his exposure to the risky assets ($\varpi_\infty^* \downarrow 0$) and investing the proceeds in the riskless asset. If the perception of the jump risk exceeds the critical value $\tilde{\lambda}$ given in (50), then the investor should go even further and start short-selling the risky assets. Because the jump risk affects all the assets, the perception of an increase in the intensity of the jumps leads the investor to dump all the risky assets indiscriminately.

4.2. *Response to jumps of different magnitudes.* If we now concentrate on the effect of an increase in the jump size instead of the jump intensity, then

$$\frac{\partial \varpi_\infty^*}{\partial \bar{J}} = \frac{1}{2\bar{J}^2}\left(1 - \frac{2\rho v^2 + \bar{J}\bar{R}}{\sqrt{(2v^2\rho - \bar{J}\bar{R})^2 + 8\bar{J}(\bar{R} + \bar{J}\lambda)v^2\rho}}\right) > 0$$

for $\bar{J} < 0$. The monotonicity implies that as the jump size gets closer to zero ($\bar{J} \uparrow 0$), the investor increases his holdings in the risky assets and conversely as $\bar{J} \downarrow (-1)$.

As to the sign of $\varpi_\infty^*$, we have

(52) $\qquad\qquad\qquad \varpi_\infty^* > 0 \qquad \text{for } -\bar{R}/\lambda < \bar{J} < 0,$

(53) $\qquad\qquad\qquad \varpi_\infty^* < 0 \qquad \text{for } -1 < \bar{J} < -\bar{R}/\lambda,$

as long as $\bar{R}/\lambda < 1$. If $\bar{R}/\lambda > 1$, the expected return is high enough relative to the jump intensity that the investor will always maintain a positive investment $\varpi_\infty^* > 0$ in the different assets, no matter how large the jump size (within the constraint $\bar{J} > -1$, of course).

4.3. *Sensitivity to risk aversion.* Here we consider the effect of the CRRA coefficient $\gamma$ on $\varpi_\infty^*$. For a CRRA investor, the first-order condition of (36) is given by

(54) $\qquad -\bar{R} + \gamma \varpi_\infty v^2 \rho - \lambda \int_0^1 \bar{J}z(1 + \varpi_\infty \bar{J}z)^{-\gamma} \nu(dz) = 0;$

then, making use of the implicit function theorem, we get

(55) $\quad \dfrac{\partial \varpi_\infty^*}{\partial \gamma} = -\dfrac{\varpi_\infty^* v^2 \rho + \lambda \int_0^1 \bar{J}z(1 + \varpi_\infty^* \bar{J}z)^{-\gamma} \ln(1 + \varpi_\infty^* \bar{J}z)\nu(dz)}{\gamma(v^2\rho + \lambda \int_0^1 \bar{J}^2 z^2 (1 + \varpi_\infty^* \bar{J}z)^{-\gamma-1}\nu(dz))}.$

The denominator is always positive but the numerator could be negative, zero or positive depending on the sign of $\varpi_\infty^*$. That is,

(56) $\qquad\qquad\qquad \dfrac{\partial \varpi_\infty^*}{\partial \gamma} = \begin{cases} > 0, & \text{if } \varpi_\infty^* < 0, \\ = 0, & \text{if } \varpi_\infty^* = 0, \\ < 0, & \text{if } \varpi_\infty^* > 0. \end{cases}$

This, in turn, implies that the higher the CRRA coefficient of an investor is, the smaller will be his $\varpi_\infty^*$ in absolute value. In the limit where $\gamma$ increases to $\infty$, $|\varpi_\infty^*|$ decreases to zero.

**5. Optimal portfolio in a multisector economy with multiple jump risks.**
We now generalize our previous results by studying the more realistic portfolio selection problem in an economy composed of $m$ sectors (or regions of the world), each containing $k$ firms (or countries). We also generalize the model to allow for multiple sources of jump risk. The total number of assets available to the investor is $n = mk$.

To understand the qualitative implications of the solution for portfolio allocation in the presence of jumps, one would presumably be primarily interested in the situation where $m$ is fixed and $k$ goes to infinity with $n$. But we provide the full solution, including the special case where all the assets are fundamentally different, that is, $k = 1$ and $m = n$. This provides a fair amount of generality to capture different situations, with different numbers of sectors.

5.1. *Sector jumps.* As in the simpler one-factor model of Section 3, we will decompose the space of returns as $\mathbb{R}^n = \bar{V} \oplus V^\perp$ where $\bar{V}$ is the span of the jump vector(s) and $V^\perp$ is the orthogonal hyperplane. With $m$ sectors, it is now natural to assume that there is potentially a separate jump term per sector, which can affect not only its own sector but also the other sectors, and do so to different degrees.

We are therefore replacing the dynamics of the risky asset returns in (2) with

$$(57) \quad \frac{dS_{i,t}}{S_{i,t-}} = (r + R_i)\, dt + \sum_{j=1}^{n} \sigma_{i,j}\, dW_{j,t} + \sum_{l=1}^{m} J_{i,l}\, dY_{l,t}, \qquad i = 1, \ldots, n,$$

where each $Y_{l,t}$ is a Lévy pure jump process with Lévy measure $\lambda_l \nu_l(dz)$, $l = 1, \ldots, m$. We assume that the Lévy pure jump processes and the Brownian motions are mutually independent and that the sum of the jumps has support on $(-1, \infty)$ to guarantee the positivity (or limited liability) of $S_i$. The constant $J_{i,l}$ is asset $i$'s scaling of sector $l$'s jump. Let $\mathbf{J}_l = [J_{1,l}, \ldots, J_{n,l}]'$.

Let us define the vector $\mathbf{1}_l$ as the $n$-vector with 1's placed in the rows corresponding to the $l$-block and zeros everywhere else:

$$(58) \quad \mathbf{1}_l = [0, \ldots, 0, \underbrace{1, \ldots, 1}_{\text{sector } l}, 0, \ldots, 0]',$$

where the first 1 is located in the $k(l-1) + 1$ coordinate, and the last one in the $kl + 1$ coordinate.

We assume that each jump vector is of the form

$$(59) \quad \mathbf{J}_l = \sum_{s=1}^{m} j_{s,l} \mathbf{1}_s = [\underbrace{j_{1,l}, \ldots, j_{1,l}}_{\text{sector } 1}, \underbrace{j_{2,l}, \ldots, j_{2,l}}_{\text{sector } 2}, \ldots, \underbrace{j_{m,l}, \ldots, j_{m,l}}_{\text{sector } m}]',$$

meaning that firms within a given sector have the same response to the arrival of a jump in a given sector, that is, to a jump in $Y_l$, but the proportional response $j_{s,l}$ of firms of different sectors to the arrival of a jump can be different, and also the proportional response of firms in a given sector to the arrival of jumps of different sectors can be different. Since some of the $j_{s,l}$ coefficients can be zero, jumps in one-sector can affect only this sector or some or all of the sectors.

So the jump vectors are unrestricted linear combinations of the $m$ vectors $\{\mathbf{1}_l\}_{l=1,\ldots,m}$. Corresponding to the above structure, we now define the orthogonal decomposition $\mathbb{R}^n = \bar{V} \oplus V^\perp$ where $\bar{V}$ is the $m$-dimensional span of the vectors $\{\mathbf{1}_l\}_{l=1,\ldots,m}$ that contains the jump vectors, and $V^\perp$ is the $(n-m)$-dimensional orthogonal hyperplane. By construction, $V^\perp$ contains no jump risk, while $\bar{V}$ contains all the jump risk. As in Section 3.1, to obtain meaningful results, we need to make the orthogonal decomposition assumption to be able to write $\Sigma = \bar{\Sigma} + \Sigma^\perp$. With this assumption, the picture is similar to Figure 1 except that $\bar{V}$ is now of dimension $m$.

5.2. *Multifactor Brownian risk.* To capture the notion of sectors, we specify a block-structure for the variance-covariance matrix of returns,

$$
(60) \quad \underset{n \times n}{\Sigma} = \begin{pmatrix} \Sigma_{1,1} & \Sigma_{1,2} & \cdots \\ \Sigma_{2,1} & \ddots & \Sigma_{2,m} \\ \cdots & \Sigma_{m,m-1} & \Sigma_{m,m} \end{pmatrix}
$$

with within-sector blocks

$$
(61) \quad \underset{k \times k}{\Sigma_{l,l}} = v_l^2 \begin{pmatrix} 1 & \rho_{l,l} & \cdots \\ \rho_{l,l} & \ddots & \rho_{l,l} \\ \cdots & \rho_{l,l} & 1 \end{pmatrix}
$$

and across-sector blocks

$$
(62) \quad \underset{k \times k}{\Sigma_{l,s}} = v_l v_s \begin{pmatrix} \rho_{l,s} & \rho_{l,s} & \cdots \\ \rho_{l,s} & \ddots & \rho_{l,s} \\ \cdots & \rho_{l,s} & \rho_{l,s} \end{pmatrix},
$$

where $1 > \rho_{l,l} > \rho_{l,s}$ and $\rho_{l,l} \geq -1/(k-1)$ and $\rho_{l,s} \geq -1/(n-1)$. In an asset pricing framework, this corresponds to a multifactor model for the returns process with $m$ common Brownian factors and $n$ idiosyncratic Brownian shocks.

Some structure on $\Sigma$ is needed, because without an explicit decomposition of $\Sigma$ there is no explicit optimal portfolio policy except to say that it is in the form $(\Sigma^\perp)^{-1}\mathbf{R}^\perp$ but we would not know what $\Sigma^\perp$ is. The block-structure of $\Sigma$ means that we have "within-sector homogeneity" of the assets, instead of the complete homogeneity in the model of Section 3. Since the number of sectors can be arbitrarily large, this is not a serious restriction and one that would have to be assumed anyway when one estimates the model lest the number of parameters becomes unwieldy.

With this structure on $\Sigma$, it turns out that we can find exactly which part of the $\Sigma$ matrix lines up with the jumps (in $\bar{V}$ defined by the specification of the $m$ jump terms), and which part is orthogonal to them (in $V^\perp$). That decomposition of $\Sigma$ is given by

$$(63) \qquad \Sigma = \underbrace{\sum_{l,s=1}^{m} \frac{\kappa_{l,s}}{k} \mathbf{1}_l \mathbf{1}'_s}_{=\bar{\Sigma}} + \Sigma^\perp,$$

where

$$(64) \qquad \kappa_{l,s} = \begin{cases} v_l^2(1 + (k-1)\rho_{l,l}), & \text{if } l = s, \\ k v_l v_s \rho_{l,s}, & \text{if } l \neq s. \end{cases}$$

$\mathbf{1}_l \mathbf{1}'_s$ is an $n \times n$ matrix with a $k \times k$ matrix of 1's placed in the $(l,s)$-block and zeros everywhere else.

The part of the $\Sigma$ matrix orthogonal to the jumps is given by the block-diagonal matrix

$$(65) \qquad \underset{n \times n}{\Sigma^\perp} = \begin{pmatrix} \Sigma^\perp_{1,1} & 0 & \cdots \\ 0 & \ddots & 0 \\ \cdots & 0 & \Sigma^\perp_{m,m} \end{pmatrix}$$

with blocks

$$(66) \qquad \underset{k \times k}{\Sigma^\perp_{l,l}} = v_l^2(1 - \rho_{l,l}) \begin{pmatrix} (k-1)/k & -1/k & \cdots \\ -1/k & \ddots & -1/k \\ \cdots & -1/k & (k-1)/k \end{pmatrix}.$$

5.3. *Optimal portfolio separation in a multisector economy.* We decompose the vector of expected excess returns on the same basis as above,

$$(67) \qquad \mathbf{R} = \sum_{l=1}^{m} r_l \mathbf{1}_l + \mathbf{R}^\perp = \bar{\mathbf{R}} + \mathbf{R}^\perp,$$

where $\mathbf{R}^\perp$ is orthogonal to each $\mathbf{1}_l$ and has the form

$$\mathbf{R}^\perp = [\mathbf{R}_1^{\perp\prime}, \ldots, \mathbf{R}_m^{\perp\prime}]'.$$

Each of the $k$-vectors $\mathbf{R}_l^\perp$ is orthogonal to the $k$-vector of 1's.

We will be looking for a vector of optimal portfolio weights in the form

$$(68) \qquad \omega = \sum_{l=1}^{m} \bar{\omega}_l \mathbf{1}_l + \omega^\perp = \bar{\omega} + \omega^\perp.$$

By the orthogonal decomposition assumption, $\Sigma^\perp$ is orthogonal to the space generated by $\{\mathbf{1}_l\}_{l=1,\ldots,m}$. Thus, the CRRA investor's minimization problem again separates as

$$(69) \quad \begin{cases} \omega^{\perp*} = \arg\min_{\omega^\perp}\{g^\perp(\omega^\perp)\}, \\ \bar{\omega}^* = \arg\min_{\bar{\omega}}\{\bar{g}(\bar{\omega})\}, \end{cases}$$

where

$$(70) \quad g^\perp(\omega^\perp) = -\omega^{\perp\prime}\mathbf{R}^\perp + \tfrac{1}{2}\gamma\omega^{\perp\prime}\Sigma^\perp\omega^\perp$$

and

$$(71) \quad \begin{aligned} \bar{g}(\bar{\omega}) &= \bar{\omega}'\bar{\mathbf{R}} + \frac{\gamma}{2}\bar{\omega}'\bar{\Sigma}\bar{\omega} + \sum_{l=1}^m \lambda_l\psi_l(\bar{\omega}'\mathbf{J}_l) \\ &= -k\sum_{l=1}^m \bar{\omega}_l r_l + \frac{\gamma}{2}k\sum_{l,s=1}^m \kappa_{l,s}\bar{\omega}_l\bar{\omega}_s + \sum_{l=1}^m \lambda_l\psi_l\left(k\sum_{s=1}^m \bar{\omega}_s j_{s,l}\right) \end{aligned}$$

with

$$(72) \quad \psi_l(x) = -\frac{1}{(1-\gamma)}\int [(1+xz)^{1-\gamma}-1]\nu_l(dz).$$

The key consequence of this separation is that in the space $V^\perp$ where only Brownian risk occurs, we again get back the Merton no-jump solution (with $\Sigma^\perp$ and $\mathbf{R}^\perp$ instead of the full $\Sigma$ and $\mathbf{R}$). Indeed, the first-order condition for minimizing $g^\perp(\omega^\perp)$ is

$$(73) \quad \Sigma^\perp \omega^{\perp*} = \frac{1}{\gamma}\mathbf{R}^\perp,$$

which, by the block diagonal form of $\Sigma^\perp$ and by construction $\mathbf{1}_l'\omega^\perp = 0$ for all $l$, leads to the optimal solution $\omega^{\perp*} = [\omega_1^{\perp*\prime},\ldots,\omega_m^{\perp*\prime}]'$ with

$$(74) \quad \omega_l^{\perp*} = \frac{1}{\gamma v_l^2(1-\rho_{l,l})}\mathbf{R}_l^\perp$$

for $l=1,\ldots,m$. Note that in general, the null space of $\Sigma^\perp$ is $\bar{V}$ and so (73) has a unique solution in $V^\perp$.

In the space $\bar{V}$, both jump and Brownian risks coexist. Minimizing $\bar{g}(\bar{\omega})$ leads to the analogue of what happens with one-sector, as in Section 3, but in dimension $m$. Similarly its solution has a limit as $k$ goes to infinity with $n$ (the number of sectors $m$ being fixed). With the change of variable $\varpi_n = k\bar{\omega}$ we see that

$$\varpi_n^* = \arg\min_{\{\varpi_n\}}\left\{-\underbrace{\sum_{l=1}^m \varpi_{nl} r_l}_{\text{Return contribution (part in } \bar{V})}\right.$$

$$
\text{(75)} \qquad + \underbrace{\frac{\gamma}{2}\sum_{l,s=1}^{m}\frac{\kappa_{l,s}}{k}\varpi_{nl}\varpi_{ns}}_{\text{Brownian risk contribution (part in }\bar{V}\text{)}}
$$

$$
+ \underbrace{\sum_{l=1}^{m}\lambda_l\psi_l\left(\sum_{s=1}^{m}\varpi_{ns}j_{s,l}\right)}_{\text{Jump risk contribution (all in }\bar{V}\text{)}} \Bigg\}
$$

which, compared to (36), is an $m$-dimensional minimization problem, instead of a one-dimensional one. The convexity of the objective function implies the existence of a unique minimizer. Letting $k \to \infty$, we see that $\varpi_n^*$ converges to a finite limit, $\varpi_\infty^*$, with the limit given by the equation above with $\kappa_{l,s}/k$ replaced by $v_l v_s \rho_{l,s}$. As in the one-factor case, we can determine below $\varpi_n^*$ in closed form under some specific jump distributions.

The wealth process $X_t^*$ of the optimizing investor will have geometric Lévy dynamics. Under the natural condition that $\mathbf{R}_l^{\perp\prime}\mathbf{R}_l^\perp = O(k)$, as $k$, the number of stocks per sector, increases, the optimal portfolio can achieve expected gains at the expense of variance which both grow approximately linearly with $k$, while keeping the exposure to the sector jumps bounded. This result is achieved by the investor apportioning an increasing fraction of assets in the subspace $V^\perp$ orthogonal to the vectors $\mathbf{1}_l$. And it follows that a generalization of the separation theorem of Section 3.4 will also hold in this situation.

5.4. *Three-fund separation.* More specifically, consider for simplicity a model in which the common jumps are identically distributed, that is, $\nu_l(dz) = \nu(dz)$ and $\mathbf{J}_l = \mathbf{J}$ for all $l = 1, \ldots, m$. The objective is to minimize

$$
\text{(76)} \qquad g(\omega) = -\omega'\mathbf{R} + \frac{\gamma}{2}\omega'\Sigma\omega + \lambda m \psi(\omega'\mathbf{J})
$$

with

$$
\text{(77)} \qquad \psi(\omega'\mathbf{J}) = -\frac{1}{(1-\gamma)}\int[(1+\omega'\mathbf{J}z)^{1-\gamma} - 1]\nu(dz).
$$

The first-order conditions for $\omega$ are given by

$$
\text{(78)} \qquad -\mathbf{R} + \gamma\Sigma\omega + m\lambda\mathbf{J}\dot{\psi}(\omega'\mathbf{J}) = 0.
$$

Define the scalar $y = \omega'\mathbf{J}$. After multiplying (78) by $\frac{1}{\gamma}\mathbf{J}'\Sigma^{-1}$ we have that $y$ must satisfy

$$
\text{(79)} \qquad -\frac{1}{\gamma}\mathbf{J}'\Sigma^{-1}\mathbf{R} + y + \frac{m\lambda}{\gamma}\mathbf{J}'\Sigma^{-1}\mathbf{J}\dot{\psi}(y) = 0;
$$

then use (78) to solve for $\omega$:

$$\omega = \frac{1}{\gamma}\Sigma^{-1}\mathbf{R} - \frac{m\lambda}{\gamma}\Sigma^{-1}\mathbf{J}\dot{\psi}(y) \tag{80}$$

or

$$\begin{aligned}\omega &= \frac{1}{\gamma}\Sigma^{-1}\mathbf{R} + \Sigma^{-1}\mathbf{J}\left(\frac{-1/\gamma\mathbf{J}'\Sigma^{-1}\mathbf{R}+y}{\mathbf{J}'\Sigma^{-1}\mathbf{J}}\right)\\ &= \frac{1}{\gamma}\underbrace{\Sigma^{-1}\left(\mathbf{R} - \frac{\mathbf{J}'\Sigma^{-1}\mathbf{R}}{\mathbf{J}'\Sigma^{-1}\mathbf{J}}\mathbf{J}\right)}_{\delta_1} + y\underbrace{\frac{1}{\mathbf{J}'\Sigma^{-1}\mathbf{J}}\Sigma^{-1}\mathbf{J}}_{\delta_2}.\end{aligned} \tag{81}$$

It is clear from the analysis above that three-fund separation continues to hold in the multisector case. The three funds are:

- a mean-variance fund which holds proportions $\delta_1 = \Sigma^{-1}(\mathbf{R} - \frac{\mathbf{J}'\Sigma^{-1}\mathbf{R}}{\mathbf{J}'\Sigma^{-1}\mathbf{J}}\mathbf{J})$ of the $n$ risky assets,
- a hedging fund holding a portfolio $\delta_2 = \frac{1}{\mathbf{J}'\Sigma^{-1}\mathbf{J}}\Sigma^{-1}\mathbf{J}$ of the $n$ risky assets; how much of this fund the investor is willing to take is controlled by the quantity $y$,
- a risk-free fund holding the riskless asset only.

When the common jumps of the $m$-region model are not identical, the three-fund separation property generalizes to an $(m+2)$-fund separation property.

5.5. *Examples of fully closed-form portfolio weights.* For the $m$-region model in which the common jumps are identically distributed, a CRRA investor with CRRA coefficient $\gamma$ will optimize the following objective function:

$$\bar{g}_n(\varpi/k) = -\sum_{l=1}^{m}\varpi_{nl}\bar{r}_l + \frac{\gamma}{2}\sum_{l=1}^{m}\sum_{s=1}^{m}\frac{\kappa_{l,s}}{k}\varpi_{nl}\varpi_{ns} \\ + m\lambda\psi\left(\sum_{s=1}^{m}\varpi_{ns}j_s\right), \tag{82}$$

where

$$\psi(y) = -\frac{1}{(1-\gamma)}\int[(1+yz)^{1-\gamma}-1]\nu(dz). \tag{83}$$

The optimal solution must satisfy the solvency constraints

$$Z_{\sup}\sum_{s=1}^{m}\varpi_{ns}j_s > -1 \quad \text{and} \quad Z_{\inf}\sum_{s=1}^{m}\varpi_{ns}j_s > -1, \tag{84}$$

where $Z_{\sup}$ and $Z_{\inf}$ are the upper and lower bounds of the support of $\nu(dz)$. Otherwise, there is a positive probability of wealth $X_t$ becoming negative, which is inadmissible in the power utility case.

The first-order conditions for $\varpi_n$ are given by

$$(85) \quad -r_l + \gamma \sum_{s=1}^{m} \frac{\kappa_{l,s}}{k}\varpi_{ns} + j_l m\lambda \dot\psi\left(\sum_{s=1}^{m}\varpi_{ns}j_s\right) = 0$$

for $l = 1,\ldots,m$. These first-order conditions form a system of $m$ equations which admit a unique solution $\varpi_n$ satisfying the solvency constraints. Let $\mathbf{j} = [j_1,\ldots,j_m]'$, $\mathbf{r} = [r_1,\ldots,r_m]'$, $\varpi_n = [\varpi_{n1},\ldots,\varpi_{nm}]'$,

$$(86) \quad \mathbf{K} = \begin{pmatrix} \kappa_{1,1} & \cdots & \kappa_{1,m} \\ \vdots & \ddots & \vdots \\ \kappa_{m,1} & \cdots & \kappa_{m,m} \end{pmatrix}$$

and $y = \sum_{s=1}^{m}\varpi_{ns}j_s = \mathbf{j}'\varpi_n$. Then (85) becomes

$$(87) \quad -\mathbf{r} + \frac{\gamma}{k}\mathbf{K}\varpi_n + m\lambda\mathbf{j}\dot\psi(y) = \mathbf{0};$$

multiplying (87) by $\frac{k}{\gamma}\mathbf{j}'\mathbf{K}^{-1}$ we get that $y$ must satisfy

$$(88) \quad -\frac{k}{\gamma}\mathbf{j}'\mathbf{K}^{-1}\mathbf{r} + y + \frac{mk\lambda}{\gamma}\mathbf{j}'\mathbf{K}^{-1}\mathbf{j}\dot\psi(y) = \mathbf{0}.$$

Once we know $\psi$ we can solve (88) for $y$ and then use (87) to solve for $\varpi_n$:

$$(89) \quad \varpi_n = \frac{k}{\gamma}\mathbf{K}^{-1}\mathbf{r} - \frac{mk\lambda}{\gamma}\mathbf{K}^{-1}\mathbf{j}\dot\psi(y)$$

or

$$(90) \quad \varpi_n = \frac{k}{\gamma}\mathbf{K}^{-1}\mathbf{r} + \mathbf{K}^{-1}\mathbf{j}\left(\frac{-k/\gamma \mathbf{j}'\mathbf{K}^{-1}\mathbf{r} + y}{\mathbf{j}'\mathbf{K}^{-1}\mathbf{j}}\right).$$

We again allow the Lévy measure to generate asymmetric positive and negative jumps. So, we consider a Lévy measure such that

$$(91) \quad m\lambda\nu(dz) = \begin{cases} \lambda_+\, dz/z, & \text{if } z \in (0,1], \\ -\lambda_-\, dz/z, & \text{if } z \in [-1,0) \end{cases}$$

with $\lambda_+ > 0$ and $\lambda_- > 0$. For a reasonable value of $\gamma$, say $\gamma = 2$, then $m\lambda\psi(y) = -\lambda_+\log(1+y) - \lambda_-\log(1-y)$ and (88) becomes a cubic equation in $y$:

$$(92) \quad f'(y) = -A + y + B[-\lambda_+(1+y)^{-1} + \lambda_-(1-y)^{-1}] = 0$$

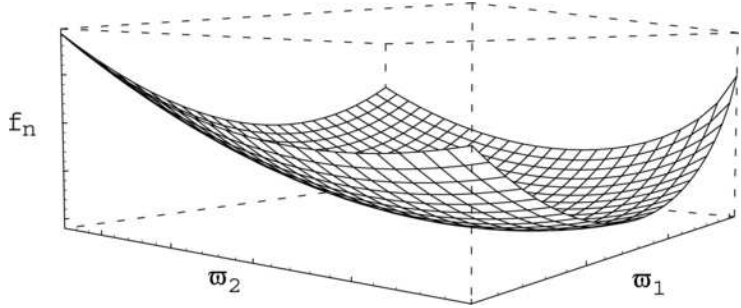

Fig. 4. *Bivariate objective function in a two-sector economy.*

with solvency constraint $|y| < 1$, where $A = \frac{k}{2}\mathbf{j}'\mathbf{K}^{-1}\mathbf{r}$ and $B = \frac{k}{2}\mathbf{j}'\mathbf{K}^{-1}\mathbf{j}$. The optimal solution to (88) under the solvency constraint is solvable fully in closed form:

$$
\begin{aligned}
y = \frac{A}{3} + \sqrt{\frac{4}{3}\left(1 + \frac{A^2}{3} + B(\lambda_+ + \lambda_-)\right)} \\
\times \cos\left(\frac{1}{3}\arccos\left(\left(\frac{A}{3}\left(-1 + \frac{A^2}{9} + \frac{3B(\lambda_- - \lambda_+)}{2A}\right.\right.\right.\right. \\
\left.\left.\left.\left. + \frac{B(\lambda_+ + \lambda_-)}{2}\right)\right)\right. \\
\left.\left. \times \left(\left(\sqrt{\frac{1}{3}\left(1 + \frac{A^2}{3} + B(\lambda_+ + \lambda_-)\right)}\right)^3\right)^{-1}\right)\right. \\
\left. + \frac{4}{3}\pi\right).
\end{aligned}
\tag{93}
$$

Multiple sectors provide the investors with additional options. For instance, if $\mathbf{r} = -(\lambda_+ - \lambda_-)\mathbf{j}$, then the investor has no exposure to jump risk. The reason is that jumps in one-sector are used to offset jump risk in another sector. Figure 4 plots the objective function, $f_n(\varpi) = \bar{g}(\varpi/k)$ that we obtain in a two-sector economy, that is the function (82) with $m = 2$.

**6. Other utility functions.** Before concluding, we briefly outline how comparable results can be derived for the other two main families of preferences, namely exponential and log utilities.

6.1. *Exponential utility.* As seen above, for an investor with CRRA utility function the optimal wealth process $X_t^*$ achieved by picking the constant portfolio fractions $\omega^*$ is itself a one-dimensional geometric Lévy process whose characteristic triple is $(X_{t-}^*(\omega^{*\prime}\mathbf{R} + r - K), X_{t-}^{*2}\omega^{*\prime}\sigma\sigma'\omega^*, f)$ where $f(dy) = \nu(dz)$ with $y = X^*_{t-}\omega'\mathbf{J}z$. In this case, the investor keeps constant fractions of wealth in each risky asset, and the constant remaining fraction $1 - \sum_j \omega_j$ in the riskless asset.

Now, consider an investor with exponential utility, $U(C) = -\frac{1}{q}\exp(-qC)$ with CARA coefficient $q > 0$. We can look for a solution to (8) in the form

$$L(x) = -\frac{K}{q}e^{-rqx}, \tag{94}$$

so that

$$\frac{\partial L(x)}{\partial x} = -rqL(x), \qquad \frac{\partial^2 L(x)}{\partial x^2} = r^2q^2L(x). \tag{95}$$

Then (8) reduces to

$$\begin{aligned}
0 = \max_{\{C_t,\omega_t\}} \Big\{ & U(C_t) - \beta L(X_t) - rqL(X_t)(X_tr + X_t\omega'_t\mathbf{R} - C_t) \\
& + \frac{1}{2}r^2q^2L(X_t)X_t^2\omega'_t\Sigma\omega_t \\
& + \lambda\int[e^{-rqX_t\omega'_t\mathbf{J}z}L(X_t) - L(X_t)]\nu(dz) \Big\},
\end{aligned} \tag{96}$$

that is,

$$\begin{aligned}
0 = \min_{\{C_t,\omega_t\}} \Big\{ & \frac{U(C_t)}{rqL(X_t)} - \frac{\beta}{qr} - (X_tr + X_t\omega'_t\mathbf{R} - C_t) \\
& + \frac{1}{2}rqX_t^2\omega'_t\Sigma\omega_t + \frac{\lambda}{rq}\int[e^{-qX_t\omega'_t\mathbf{J}z} - 1]\nu(dz) \Big\}
\end{aligned} \tag{97}$$

after division by $qL(X_t)$ [note that max becomes min as a result of $qL(X_t) < 0$].

The optimal policy of $\omega = X_t\omega_t$ is given by the objective function

$$\min_{\{\omega\}} \Big( -\omega'\mathbf{R} + \frac{1}{2}rq\omega'\Sigma\omega + \frac{\lambda}{rq}\int[e^{-rq\omega'\mathbf{J}z} - 1]\nu(dz) \Big) \tag{98}$$

and the optimal consumption choice is therefore

$$C_t^* = rX_t - \frac{1}{q}\log(rK). \tag{99}$$

Finally, we evaluate (97) at $C^*$ and $\omega^*$ to identify $K$,

$$K = \frac{1}{r}\exp\Big(1 - \frac{\beta}{r} - q\omega^{*\prime}\mathbf{R} + \frac{1}{2}rq^2\omega^{*\prime}\Sigma\omega^* + \frac{\lambda}{r}\int[e^{-q\omega^{*\prime}\mathbf{J}z} - 1]\nu(dz)\Big). \tag{100}$$

For an investor with CARA utility function the optimal wealth process $X_t^*$ achieved by picking constant portfolio amounts is now a one dimensional arithmetic Lévy process. The investor keeps constant amounts of wealth in each risky asset, and the remaining amount in the riskless asset.

6.2. *Log utility.* Finally, consider an investor with log utility, $U(x) = \log(x)$. We can look for a solution to (8) in the form

(101) $$L(x) = K_1^{-1} \log(x) + K_2,$$

where $K_1$ and $K_2$ are constant, so that

(102) $$\frac{\partial L(x)}{\partial x} = K_1 x^{-1}, \qquad \frac{\partial^2 L(x)}{\partial x^2} = -K_1 x^{-2}.$$

Then (8) reduces to

(103) $$\begin{aligned}0 = \max_{\{C_t, \omega_t\}} \Big\{ &\log(C_t) - \beta K_1^{-1} \log(X_t) - \beta K_2 \\ &+ K_1^{-1} X_t^{-1}(rX_t + \omega_t' \mathbf{R} X_t - C_t) \\ &- \frac{1}{2} K_1^{-1} \omega_t' \Sigma \omega_t + \lambda K_1^{-1} \int \log(1 + \omega_t' \mathbf{J} z) \nu(dz) \Big\},\end{aligned}$$

that is,

(104) $$\begin{aligned}0 = \min_{\{C_t, \omega_t\}} \Big\{ &-\log(C_t) + \beta K_1^{-1} \log(X_t) + \beta K_2 \\ &- K_1^{-1} r - K_1^{-1} \omega_t' \mathbf{R} + K_1^{-1} X_t^{-1} C_t \\ &+ K_1^{-1} \frac{1}{2} \omega_t' \Sigma \omega_t - \lambda K_1^{-1} \int \log(1 + \omega_t' \mathbf{J} z) \nu(dz) \Big\}.\end{aligned}$$

The optimal policy of $\omega_t$ is given by the objective function,

(105) $$\min_{\{\omega_t\}} \Big( -\omega_t' \mathbf{R} + \frac{1}{2} \omega_t' \Sigma \omega_t - \lambda \int \log(1 + \omega_t' \mathbf{J} z) \nu(dz) \Big)$$

and the optimal consumption choice is therefore

(106) $$C_t^* = K_1 X_t.$$

To identify $K_1$ and $K_2$, we evaluate (104) at $C^*$ and $\omega^*$,

$$K_2 = \frac{1}{\beta} \Big\{ \log(\beta) + \frac{r}{\beta} + \frac{1}{\beta} \omega^{*\prime} \mathbf{R} - 1 - \frac{1}{2\beta} \omega^{*\prime} \Sigma \omega^* + \frac{\lambda}{\beta} \int \log(1 + \omega^{*\prime} \mathbf{J} z) \nu(dz) \Big\},$$

$$K_1 = \beta.$$

**7. Conclusions.** We have proposed a new approach to characterize in closed form the portfolio selection problem for an investor concerned with the possibility of jumps in asset returns, and who seeks to control this risk by diversification or other means. We extended the standard multiasset geometric Brownian motion models to exponential Lévy models through the inclusion of correlation effects due to jumps.

Our key decomposition is one of the space of returns into a space containing both jump and Brownian risks, and one containing only Brownian risk. The model is quite flexible and can capture many different situations, from jumps of different signs to jumps of a single sign, multiple jumps that can affect one or some of the sectors; it is as flexible as it can get while achieving a closed-form solution, which is critical if one is to solve a portfolio choice problem in practice with a large number of assets.

With explicit optimal portfolios policies, one can address important practical questions. How exactly does increasing the number of available assets improve the investor's exposure to both diffusive and jump risk? How does the portfolio of an investor who fears jumps differ from the portfolio of one who does not? Is there a simple form for the optimal portfolio which is achieved asymptotically as the number of assets grows to infinity? And so on.

Finally, we discuss further possible extensions and generalizations of this work:

1. One could consider modeling the spectral decomposition of the $\Sigma$ matrix directly, instead of parametrizing the matrix itself and then determining its decomposition across the two spaces $\bar{V}$ and $V^\perp$. As an empirical strategy, one could determine the number of sectors through factor analysis or similar techniques, in order to determine in a data-driven manner the shape of the $\Sigma$ matrix.
2. Stochastic volatility of the type considered in Liu, Longstaff and Pan (2003) requires solving our nonlinear equations for weight vectors stepwise in time, in parallel with ordinary differential equations (which themselves depend on the current portfolio weights). This does not appear doable in closed form.
3. Portfolio restrictions such as short-selling constraints are relevant in practice but, when generic constraints are imposed on the optimal portfolio, we cannot expect the dimensional reduction to be preserved or our conclusions to hold. However, a utility function such as power utility which becomes $-\infty$ for wealth below a finite threshold, sometimes automatically implies certain constraints; it appears that in this case, much of our analysis remains intact.
4. Finally, we would like to be able to capture a more subtle form of contagion, in the form not just of simultaneous jumps within or across sectors, as we are currently able to model, but rather in the form of a jump in one-sector causing an increase in the likelihood that a different jump will occur in another sector. To capture this effect, self-exciting jump processes seem a promising approach which we intend to investigate in future work.

**Acknowledgments.** We are very grateful to Jun Liu, Francis Longstaff and Raman Uppal for comments. We are also grateful to seminar participants at LBS, LSE, Stanford, Tilburg, UIC, UIUC, Wharton, the European Science Foundation Conference on Advanced Mathematical Methods for Finance and the Econometric Society Australasian Meeting (E. J. Hannan Lecture).

## REFERENCES

Aase, K. K. (1984). Optimum portfolio diversification in a general continuous-time model. *Stochastic Process. Appl.* **18** 81–98. MR757349

Ang, A. and Bekaert, G. (2002). International asset allocation with regime shifts. *Review of Financial Studies* **15** 1137–1187.

Ang, A. and Chen, J. (2002). Asymmetric correlations of equity portfolios. *Journal of Financial Economics* **63** 443–494.

Bae, K.-H., Karolyi, G. A. and Stulz, R. M. (2003). A new approach to measuring financial contagion. *Review of Financial Studies* **16** 717–763.

Choulli, T. and Hurd, T. R. (2001). The role of Hellinger processes in mathematical finance. *Entropy* **3** 150–161.

Cvitanić, J., Polimenis, V. and Zapatero, F. (2008). Optimal portfolio allocation with higher moments. *Annals of Finance* **4** 1–28.

Das, S. and Uppal, R. (2004). Systemic risk and international portfolio choice. *Journal of Finance* **59** 2809–2834.

Emmer, S. and Klüppelberg, C. (2004). Optimal portfolios when stock prices follow an exponential Lévy process. *Finance Stoch.* **8** 17–44. MR2022977

Grauer, R. and Hakansson, N. (1987). Gains from international diversification: 1968–1985 returns on portfolios of stocks and bonds. *Journal of Finance* **42** 721–739.

Han, S. and Rachev, S. (2000). Portfolio management with stable distributions. *Math. Methods Oper. Res.* **51** 341–352. MR1761863

Hartmann, P., Straetmans, S. and de Vries, C. (2004). Asset market linkages in crisis periods. *Review of Economics and Statistics* **86** 313–326.

Jeanblanc-Picqué, M. and Pontier, M. (1990). Optimal portfolio for a small investor in a market model with discontinuous prices. *Appl. Math. Optim.* **22** 287–310. MR1068184

Kallsen, J. (2000). Optimal portfolios for exponential Lévy processes. *Math. Methods Oper. Res.* **51** 357–374. MR1778648

Liu, J., Longstaff, F. and Pan, J. (2003). Dynamic asset allocation with event risk. *J. Finance* **58** 231–259.

Longin, F. and Solnik, B. (2001). Extreme correlation of international equity markets. *J. Finance* **56** 649–676.

Madan, D. (2004). Equilibrium asset pricing with non-Gaussian factors and exponential utilities. Technical report, Univ. Maryland.

Merton, R. C. (1969). Lifetime portfolio selection under uncertainty: The continuous-time case. *Review of Economics and Statistics* **51** 247–257.

Merton, R. C. (1971). Optimum consumption and portfolio rules in a continuous-time model. *J. Econom. Theory* **3** 373–413. MR0456373

Ortobelli, S., Huber, I., Rachev, S. T. and Schwartz, E. S. (2003). Portfolio choice theory with non-Gaussian distributed returns. In *Handbook of Heavy Tailed Distributions in Finance* (S. T. Rachev, ed.) 547–594. Elsevier, Amsterdam.

Shirakawa, H. (1990). Optimal dividend and portfolio decisions with Poisson and diffusion-type return process. Technical report, Tokyo Institute of Technology.

Solnik, B. (1974). Why not diversify internationally rather than domestically? *Financial Analysts Journal* **30** 48–53.

Y. Aït-Sahalia
Department of Economics
Princeton University
and
NBER
Princeton, New Jersey 08544-1021
USA
E-mail: yacine@princeton.edu

J. Cacho-Diaz
Department of Economics
Princeton University
Princeton, New Jersey 08544-1021
USA
E-mail: jcacho@princeton.edu

T. R. Hurd
Department of Mathematics
  and Statistics
McMaster University
Hamilton, Ontario L8S 4K1
Canada
E-mail: hurdt@mcmaster.ca